\documentclass[11pt]{article}

\usepackage[margin=1in]{geometry}
\usepackage{amsmath,amssymb,amsthm}
\usepackage{array}
\usepackage{booktabs}
\usepackage{tikz}
\usepackage{graphicx}
\usepackage{hyperref}
\usepackage[margin=1in]{geometry}

\theoremstyle{plain}
\newtheorem{theorem}{Theorem}[section]
\newtheorem{proposition}[theorem]{Proposition}
\newtheorem{lemma}[theorem]{Lemma}
\newtheorem{corollary}[theorem]{Corollary}

\theoremstyle{definition}
\newtheorem{definition}[theorem]{Definition}

\theoremstyle{remark}
\newtheorem{remark}[theorem]{Remark}

\newtheorem{problem}[theorem]{Problem}

\newcommand{\F}{\mathbb F}
\newcommand{\C}{\mathcal C}
\newcommand{\wt}{\operatorname{wt}}

\newcommand{\rank}{\operatorname{rank}}

\newcommand{\eps}{\varepsilon}

\newcommand{\FF}{\mathbb{F}}
\newcommand{\CA}{\mathcal{C}_A}

\newcommand{\iot}{\iota}
\newcommand{\GL}{\mathrm{GL}}
\newcommand{\onev}{\mathbf{1}}

\newcommand{\Hull}{\operatorname{Hull}}
\newcommand{\Row}{\operatorname{Row}}

\title{Extremal diagonal-coherent binary codes: shortened Hamming cores, deep holes, and the 64 I-Ching hexagrams}

\author{Jon-Lark Kim \\ \\
Department of Mathematics \\
Institute for Mathematical and Data Sciences \\
Sogang University, \\ Seoul, Korea\\
jlkim@sogang.ac.kr}

\date{}

\begin{document}

\maketitle

\begin{abstract}
Each of the 64 hexagrams of the I-Ching, or Book of Changes, is a stack of six Yin–Yang lines and is therefore a binary vector of length 6; the eight hexagrams formed by repeating a single trigram span a doubled-trigram subspace
$\Delta$. We ask which binary
$[6,3,3]$ codes are compatible with this trigram structure. We show that a
$[6,3,3]$ code meets $\Delta$ in at most four codewords, that the four are then forced to be the doubled trigrams Earth, Fire, Wind and Lake, and that exactly two codes attain this bound: the shortened Hamming code
$\CA$ and its dual.
 Since $\CA$ has covering radius 2, its eight codewords, their radius-one balls, and its eight deep holes together reconstruct all 64 hexagrams.
In general, we call a binary code of length
$2m$ with $d \ge 3$ {\em diagonal-coherent} if it contains the half of the diagonal subspace $\{(u,u)\}$ coming from even-weight vectors $u$ --- the largest part of that subspace a single-error-correcting code can contain --- and {\em extremal} if its dimension is as large as possible.
For $m\ge2$ and $r=\lceil\log_2(2m+1)\rceil$,
we determine the number $N(m)$ of extremal diagonal-coherent codes with parameters $[\,2m,\ 2m-r,\ \ge3\,]$, show that diagonal coherence entails no loss of redundancy, and determine their minimum distances.
 We also put these codes in a normal form, determine their hulls, and show that each has covering radius $2$ with exactly
$2^r-1-2m$ deep-hole cosets. Finally we propose three open problems.
\end{abstract}


\medskip


\section{Introduction}
\label{intro}

The \emph{I-Ching} (\emph{Yijing}, \emph{Book of Changes})~\cite{Wilhelm1967} consists of two
historically distinct layers. The older layer, the \emph{Zhouyi}, is a
divination manual comprising, for each of the 64 hexagrams, a short hexagram
statement together with six line statements; on linguistic and historical
grounds it is generally dated to the Western Zhou, around the ninth century
BCE~\cite{Shaughnessy1996,Rutt1996}. The younger layer is a set of ten
commentaries known as the Ten Wings, among them the influential \emph{Xici
zhuan} (Great Commentary), composed over the Warring States and early Han
periods, roughly the fourth to second centuries BCE. It is this second layer
that supplies the cosmological vocabulary of Yin and Yang and that turned a
manual of divination into a philosophical classic~\cite{Smith2008}.

Each of the 64 hexagrams of the I-Ching is a stack of six Yin–Yang lines. Writing
0 for a broken Yin line and 1 for a solid Yang line identifies the hexagram set with the binary vector space
\[
        \F_2^6.
\]

This identification of hexagrams with binary vectors of length 6 has a documented origin.
Leibniz~\cite{Ryan1996} described his binary arithmetic to Joachim Bouvet
in a letter of 15 February 1701; Bouvet, then at the Kangxi court in Beijing,
recognised the correspondence with the Xiantian diagram and replied on 4
November 1701, enclosing a woodcut of Shao Yong's arrangement~\cite{Nielsen2003}. The letter
reached Leibniz only on 1 April 1703, and within days he submitted to the Paris
Academy his \emph{Explication de l'arithm\'etique binaire}, whose title
explicitly advertises that binary numeration ``gives the sense of the ancient
Chinese figures of Fohy(Fuxi)''~\cite{Leibniz1703}.


Even though such a nice connection was known since Leibniz, there have been only a few papers on mathematical structures of the 64 hexagrams.
In 1975, Goldenberg~\cite{Goldenberg} gave an elementary algebraic treatment of the I-Ching for non-mathematicians. Recently, Hu et al.~\cite{HuPetoukhovPetukhova2017} noticed that the system of I-Ching can be represented by the tables with dyadic groups of 4 bigrams, 8 trigrams and 64 hexagrams, which were believed to be fundamental archetypes of the Nature. The ancient Chinese could not have known about the genetic code of protein sequences of amino acids: in particular, the genetic code is constructed on DNA molecules using 4 nitrogenous bases, 16 doublets, and 64 triplets. Hence the authors~\cite{HuPetoukhovPetukhova2017} described the usage of dyadic groups related to the hexagrams as a foundation of the bio-mathematical doctrine of the geno-logic code.
Li~\cite{YanzhongLi} introduced a novel computational perspective based on two's complement format of a negative integer to reinterpret the mathematical essence of the I Ching hexagrams and the topological mapping between the Fuxi Square and Circular Diagrams.


\medskip

In this paper we show that coding theory provides an efficient way of understanding the 64 hexagrams. More precisely, we show that at most four of the eight doubled-trigram hexagrams can lie in a binary
$[6,3,3]$ code, that the four are then forced to be
 $$\mathcal H=\langle\,101101,\;011011\,\rangle_{\mathbb F_2}=\{\,000000,\;011011,\;101101,\;110110\,\},$$
 and that exactly two codes attain this bound: the shortened Hamming
$[6,3,3]$ code $\CA$ and its dual. For each of these two codes, the eight codewords, the balls of radius one around them, and the eight deep holes together account for all $64$ hexagrams.

We then generalize this picture from hexagrams to stacks of
$2m$ lines.
We call a binary code of length
$2m$ with $d \ge 3$ {\em diagonal-coherent} if it contains the half of the diagonal subspace $\{(u,u)\}$ coming from even-weight vectors $u$ --- the largest part of that subspace a single-error-correcting code can contain --- and {\em extremal} if its dimension is as large as possible.
For $m\ge2$ and $r=\lceil\log_2(2m+1)\rceil$, we determine the number $N(m)$ of extremal diagonal-coherent binary codes with parameters $[\,2m,\ 2m-r,\ \ge3\,]$ which by definition contains $\iot(E_m)$, the image under the doubling map $\iot(u)=(u, u)$ of the even-weight code $E_m \le \F_2^m$; Requiring it costs nothing: the dimension $2m-r$ is already the largest possible for any binary code of length $2m$ with minimum distance at least $3$.

We also determine the minimum distances of these codes: they are always 3 or 4, and the value 4 occurs precisely when $m$ is a power of $2$. Even then it is exceptional --- we obtain a closed formula for the number
$N_4(m)$ of such codes (Proposition~\ref{prop:countd4}), which gives for instance
$N_4(4)=2$ out of $N(4)=10$, and
$N_4(8)=480$ out of
$N(8)=205\,920$. We give these codes in a normal form, determine their hulls, and show that each has covering radius $2$ with exactly
$2^r-1-2m$ deep-hole cosets.

\medskip

The paper consists of six sections. In Section~\ref{sec:background}, we give a structural background of I-Ching.
In Section~\ref{sec:prelim}, we begin with basic definitions of coding theory.
Section~\ref{sec-shortHamming} constructs the core
$\CA$ and proves the uniqueness statement above.
 In Section~\ref{deepholes}, we find the eight deep holes so that the 64 hexagrams can be described by $\CA$ and the deep holes. Finally, Section~\ref{sec:general} considers the generalization of Section~\ref{sec-shortHamming} by describing codes coherent with a diagonal subspace, together with three open problems.


\section{Structural background of I-Ching}
\label{sec:background}

This section collects the background needed to read the combinatorial results
below, and fixes the terminology that recurs throughout the paper.
Readers familiar with the \emph{Book of Changes} may proceed directly to
Section~\ref{sec:prelim}.

\subsection{Lines, trigrams, and hexagrams}

The atomic symbol is a \emph{line} of which there are two kinds:
a broken \emph{Yin} line and an unbroken \emph{Yang} line. Three lines stacked
vertically form a \emph{trigram}, and six form a \emph{hexagram}; lines are
counted from the bottom upward, so that position~1 is the lowest line and
position~6 the highest. There are $2^3=8$ trigrams and $2^6=64$ hexagrams, and
every hexagram decomposes uniquely into a \emph{lower} trigram (positions
$1,2,3$) and an \emph{upper} trigram (positions $4,5,6$); the 64 hexagrams are
correspondingly the $8\times 8$ ordered pairs of trigrams. Each trigram carries
a name and a canonical natural image:

\begin{center}
\begin{tabular}{>{\centering\arraybackslash}p{2.2cm}
                >{\centering\arraybackslash}p{2.2cm}
                >{\centering\arraybackslash}p{2.4cm}
                >{\centering\arraybackslash}p{4.0cm}}
\hline
Name & Image & Trigrams $u$ & Doubled-trigram hexagrams $(u, u)$\\
\hline
Qian   & Heaven   & $111$ & $111111$\\
Dui    & Lake     & $110$ & $110110$\\
Li     & Fire     & $101$ & $101101$\\
Zhen   & Thunder  & $100$ & $100100$ \\
Xun    & Wind     & $011$ & $011011$\\
Kan    & Water    & $010$ & $010010$\\
Gen    & Mountain & $001$ & $001001$\\
Kun    & Earth    & $000$ & $000000$\\
\hline
\end{tabular}
\end{center}

\noindent
Here the binary column already uses the encoding fixed in Section~\ref{intro}, namely Yin
$=0$ and Yang $=1$, with the bottom line written first. A hexagram is named by
its two trigram images; for instance Fire over Water is the hexagram called Wei Ji (``Before Completion''); note that the traditional name lists the upper trigram first, so in our encoding, which writes the lower trigram first, this is
$010101$ — Water below, Fire above.

\subsection{Divination and the notion of change}

The name \emph{Book of Changes} refers to a procedure, not merely to a text.
In the classical yarrow-stalk method --- and in the simpler three-coin method
that later displaced it --- each of the six lines is generated independently,
and the outcome distinguishes four cases rather than two: a line is Yin or
Yang, and additionally \emph{stable} or \emph{changing}. A consultation
therefore returns not one hexagram but two: the hexagram obtained by reading
the lines as drawn, and a second hexagram obtained from it by flipping exactly
the changing lines.

This is the point at which the tradition itself already performs the operation
that the present paper formalises. If $h$ denotes the first hexagram and
$S\subseteq\{1,\dots,6\}$ the set of changing lines, then the second hexagram
is
\[
  h' \;=\; h + \chi_S ,
\]
where $\chi_S$ is the characteristic vector of $S$ and addition is
componentwise modulo~2. In other words, divinatory change \emph{is} translation
in $\FF_2^6$, and the number of changing lines is the Hamming distance
$d_H(h,h')$.

\subsection{Classical relations between hexagrams}

Traditional hermeneutics singles out two operations relating one hexagram to
another; they are involutions of $\FF_2^6$ and are used repeatedly
below.

\begin{itemize}
\item \emph{Complementation} (\emph{cuo}): every line is
      exchanged for its opposite. In binary terms $h\mapsto h+111111$.
\item \emph{Overturning} (\emph{zong}): the hexagram is read
      upside down, so that position $i$ is exchanged with position $7-i$. This
      is the coordinate reversal $h_1h_2h_3h_4h_5h_6 \mapsto
      h_6h_5h_4h_3h_2h_1$. Exactly eight hexagrams are unchanged by it, namely
      the palindromic ones.
\end{itemize}

Complementation and overturning commute and generate a Klein four-group acting
on the 64 hexagrams; this action governs the pairing structure of the traditional King Wen sequence. That the dual of our
central code is realised by \emph{overturning} rather than by an arbitrary
coordinate permutation is the content of Section~\ref{subsec:uniqueness}.

\section{Preliminaries}\label{sec:prelim}

We refer to~\cite{HuffmanPless},~\cite{MacWilliamsSloane} for basic facts from Coding Theory.

Let $\F_2$ be a finite field with two elements 0 and 1.  A {\it linear code of length $n$ over $\F_2$} is an $\F_2$-linear subspace
\[
        \C\subseteq \F_2^n.
\]

The {\it Hamming weight} of a vector $x=(x_1,\ldots,x_n)\in \F_2^n$ is
$
        \wt(x)=|\{i:x_i\neq 0\}|,
$
and the {\it Hamming distance between} two vectors is
$
        d(x,y)=\wt(x-y).
$
The {\it minimum distance} of a code $\C$ is
\[
        d(\C)=\min\{\wt(c):0\neq c\in \C\}.
\]

If $\dim_{\F_2} \C=k$ and the {\em minimum Hamming distance} of $\C$ is $d$, then $\C$ has parameters $[n,k,d]$.


The \emph{covering radius} of a set $S\subseteq\FF_2^n$ is
\[
  \rho(S)=\max_{x\in\FF_2^n}\ \min_{s\in S} d(x,s);
\]
equivalently, $\rho(S)$ is the least $R$ for which the balls of radius $R$
centred at the points of $S$ cover $\FF_2^n$. Following~\cite{CohenHonkalaLitsynLobstein}, we say that
$S$ \emph{covers with radius} $R$ when $\rho(S)\le R$.

A matrix $H_\C\in \F_2^{r\times n}$ is a {\it parity-check matrix} for $\C$ if
\[
        \C=\ker H_\C=\{x\in \F_2^n:H_\C~ x^T=0\}.
\]
The dimension is then
\[
        \dim_{\F_2} \C=n-\rank_{\F_2} H_\C.
\]
The dual code is
\[
        \C^\perp=\{y\in \F_2^n:\langle x,y\rangle=0\text{ for all }x\in \C\},
\]
where
\[
        \langle x,y\rangle=\sum_{i=1}^n x_i y_i.
\]
If $\C=\ker H_\C$, then $\C^\perp$ is the row space of $H_\C$.

A linear code $\C$ is {\it self-dual} if $\C = \C^{\perp}$ and {\it self-orthogonal} if
$\C \subseteq \C^{\perp}$.
A linear code $\C$ is {\it isodual} if there exists a permutation $\sigma$ such that
$\sigma(\C) =\C^{\perp}$. It is easy to see that a self-dual code is an isodual code, but not the converse.

\section{A shortened-Hamming core with I-Ching structure} \label{sec-shortHamming}

In this section, we describe the 64 hexagrams in terms of a symmetric binary $[6, 3, 3]$ code; its deep holes are treated in Section~\ref{deepholes}.
This gives a new way to visualize the 64 hexagrams.

\subsection{A binary $[6,3,3]$ code}

We now describe a particularly symmetric binary \([6,3,3]\) code that can be used as a
central set of \(8\) hexagrams. Note that this code is equivalent to the shortened Hamming $[6,3,3]$ code~\cite{HuffmanPless} and is distance-optimal~\cite{Grassl}.

 We use the convention that a hexagram is a vector
\[
        h=(h_1,h_2,h_3,h_4,h_5,h_6)\in \mathbb F_2^6,
\]
where the first three coordinates form the lower trigram and the last three coordinates
form the upper trigram.  Thus a doubled-trigram hexagram has the form
\[
        (u,u)\in \mathbb F_2^3\oplus \mathbb F_2^3.
\]

Consider the binary linear code $\C_A\leq \mathbb F_2^6$ generated by
\[
        101101,\qquad 011011,\qquad 000111.
\]
Equivalently, $\C_A$ has generator matrix
\[
        G_A=
        \begin{pmatrix}
        1&0&1&1&0&1\\
        0&1&1&0&1&1\\
        0&0&0&1&1&1
        \end{pmatrix}.
\]

\begin{remark}
It is easy to see that $\CA^{\perp}=\langle 101101,\,011011,\,111000\rangle_{\F_2}$.
\end{remark}

\begin{proposition} \label{prop-[6,3,3]}
The code $\C_A$ is a binary linear $[6,3,3]$ code.
\end{proposition}

\begin{proof}
The three generator rows are linearly independent, so $\dim \C_A=3$.
The eight codewords are
\[
\begin{aligned}
\C_A=\{&
000000,\ 000111,\ 011011,\ 011100,\\
&101010,\ 101101,\ 110001,\ 110110
\}.
\end{aligned}
\]
The nonzero weights are $3$ or $4$.  Hence the minimum distance is
\[
        d(\C_A)=3.
\]
Therefore $\C_A$ has parameters $[6,3,3]$.
\end{proof}

Since $d(\C_A)=3$, the radius-one balls around the codewords of $\C_A$ are disjoint.
Equivalently, if
\[
        e_1=100000,\quad e_2=010000,\quad e_3=001000,\quad
        e_4=000100,\quad e_5=000010,\quad e_6=000001,
\]
then the seven sets
\[
        \C_A,\quad \C_A+e_1,\quad \C_A+e_2,\quad \C_A+e_3,\quad
        \C_A+e_4,\quad \C_A+e_5,\quad \C_A+e_6
\]
are pairwise disjoint.  Hence they cover
\[
        8\cdot 7=56
\]
of the $64$ hexagrams.

The remaining $8$ hexagrams form a single coset.  For this choice of $\C_A$, the missing
coset is
\[
        M_A=001001+ \C_A.
\]
Explicitly,
\[
\begin{aligned}
M_A=\{&
001001,\ 001110,\ 010010,\ 010101,\\
&100011,\ 100100,\ 111000,\ 111111
\}.
\end{aligned}
\]
Thus the full binary hexagram space decomposes as
\[
        \mathbb F_2^6
        =
        \C_A
        \sqcup
        \left(\bigsqcup_{i=1}^6 (\C_A+e_i)\right)
        \sqcup
        M_A.
\]

On the other hand, since
\[
        e_3+e_6=001001,
\]
we may take
\[
        M_A=\C_A+(e_3+e_6)=\C_A+001001.
\]

Thus the whole space of $64$ hexagrams has the compact coset
description
\[
\boxed{
        \mathbb F_2^6
        =
        \C_A+\{0,e_1,e_2,e_3,e_4,e_5,e_6,e_3+e_6\}.
}
\]
Equivalently,
\[
\boxed{
\mathbb F_2^6
=
\C_A
\sqcup
(\C_A+e_1)
\sqcup
\cdots
\sqcup
(\C_A+e_6)
\sqcup
(\C_A+e_3+e_6).
}
\]

\subsection{Uniqueness of the trigram-coherent core}
\label{subsec:uniqueness}

The code $\CA$ of Proposition~\ref{prop-[6,3,3]} was written down explicitly, and one may
reasonably ask whether the I-Ching features it exhibits are an artifact of that
particular choice. They are not. In this subsection we show that $\CA$ and its
dual are the \emph{only} binary $[6,3,3]$ codes that are compatible with the
trigram decomposition, in the precise sense of containing the largest possible
number of doubled-trigram hexagrams. The proof is elementary and requires no
computer search.

Throughout, a hexagram is a vector $h=(h_1,\dots,h_6)\in\FF_2^6$, the first
three coordinates forming the lower trigram and the last three the upper
trigram, and we write $\FF_2^6=\FF_2^3\oplus\FF_2^3$ accordingly. Define the
\emph{doubled-trigram subspace}
\[
  \Delta \;=\; \{\,(u,u) \;:\; u\in\FF_2^3 \,\}\;\le\;\FF_2^6 ,
  \qquad |\Delta| = 8,
\]
whose eight elements are exactly the eight hexagrams built from a single
trigram repeated twice, and the coordinate reversal
\[
  \tau(h_1h_2h_3h_4h_5h_6) \;=\; h_6h_5h_4h_3h_2h_1 ,
\]
which is the I-Ching \emph{overturning} (\emph{zong}) operation. Finally put
\[
  \mathcal H \;=\; \langle\, 101101,\; 011011 \,\rangle_{\F_2}
      \;=\; \{\,000000,\;011011,\;101101,\;110110\,\},
\]
the doubled tetrad $\{\text{Earth},\text{Fire},\text{Wind},\text{Lake}\}$.

\begin{definition}\label{def:coherent}
A binary $[6,3,3]$ code $\C\le\FF_2^6$ is \emph{trigram-coherent} if
$|\C\cap\Delta|$ is as large as possible among all binary $[6,3,3]$ codes.
\end{definition}

In the terminology of Section~\ref{sec:general}, a trigram-coherent
$[6,3,3]$ code is precisely an extremal diagonal-coherent code of length
$6$ (Definition~\ref{def:coherentgen}).

\medskip

We first bound $|\C\cap\Delta|$. Write
\[
  E=\{000,011,101,110\},\qquad O=\{100,010,001,111\}
\]
for the even- and odd-weight cosets in $\FF_2^3$, so that $E\le\FF_2^3$ is the
even-weight subspace and $O=E+111$.

\begin{lemma}\label{lem:subgroup}
Let $T\le\FF_2^3$ be a subgroup containing no vector of Hamming weight $1$.
Then $|T|\le 4$, and $|T|=4$ if and only if $T=E$.
\end{lemma}

\begin{proof}
If $\dim T=3$ then $T=\FF_2^3$ contains $100$, which is excluded. Suppose
$\dim T=2$. Then $T=\ker\varphi$ for a nonzero linear functional
$\varphi(x)=a\cdot x$ with $a\in\FF_2^3\setminus\{0\}$, and $T$ avoids the
three weight-one vectors $e_1,e_2,e_3$ if and only if $a\cdot e_i=1$ for
$i=1,2,3$, that is, if and only if $a=111$. Hence $T=\ker(x\mapsto
x_1+x_2+x_3)=E$ is the unique such subgroup of dimension $2$, and no subgroup
of larger order qualifies.
\end{proof}

\begin{proposition}\label{prop:bound}
Let $\C\le\FF_2^6$ be a binary $[6,3,3]$ code. Then
\[
  |\C\cap\Delta| \;\le\; 4 ,
\]
and equality holds if and only if $\C\cap\Delta=\mathcal H$. In particular a
trigram-coherent code contains exactly four doubled trigrams, namely the tetrad
$\{\text{Earth},\text{Fire},\text{Wind},\text{Lake}\}$.
\end{proposition}

\begin{proof}
The map $\iota:\FF_2^3\to\Delta$, $u\mapsto(u,u)$, is a group isomorphism and
satisfies $\wt(\iota(u))=2\wt(u)$. Consequently $\iota(u)$ has Hamming weight
$2$ precisely when $\wt(u)=1$.

Now $\C\cap\Delta$ is a subgroup of $\Delta$, so $T:=\iota^{-1}(\C\cap\Delta)$ is
a subgroup of $\FF_2^3$. Since $d(\C)=3$, the code $\C$ contains no nonzero
vector of weight $2$, hence $T$ contains no vector of weight $1$. By
Lemma~\ref{lem:subgroup}, $|T|\le 4$ with equality if and only if $T=E$, and
$\iota(E)=\mathcal H$.
\end{proof}

Thus the maximum in Definition~\ref{def:coherent} is $4$, and every
trigram-coherent code contains $\mathcal H$. The next theorem shows that this forces the
code completely.

\begin{theorem}\label{thm:uniqueness}
There are exactly two binary $[6,3,3]$ codes containing $\mathcal H$, namely
\[
  \CA=\langle 101101,\,011011,\,000111\rangle
  \qquad\text{and}\qquad
  \CA^{\perp}=\langle 101101,\,011011,\,111000\rangle .
\]
Equivalently, $\CA$ and $\CA^{\perp}$ are the only trigram-coherent binary
$[6,3,3]$ codes, and they satisfy $\CA\cap\CA^{\perp}=\mathcal H$.
\end{theorem}

\begin{proof}
Let $\C$ be a $[6,3,3]$ code with $\mathcal H \subseteq \C$. Since $\dim \mathcal H=2$ and
$\dim \C=3$, we have $\C=\mathcal H\sqcup(\mathcal H+c)$ for some $c=(a,b)\in\FF_2^6\setminus \mathcal H$
with $a,b\in\FF_2^3$. Writing $\mathcal H=\{(u,u):u\in E\}$, the nontrivial coset is
\[
  \mathcal H+c \;=\; \bigl\{\,(a+u,\;b+u)\;:\;u\in E\,\bigr\},
\]
and, because the nonzero words of $\mathcal H$ have weight $4$, the condition $d(\C)=3$
is equivalent to
\begin{equation}\label{eq:cosetcond}
  \wt(a+u)+\wt(b+u)\;\ge\;3 \qquad\text{for all } u\in E .
\end{equation}
Note that for $u\in E$ the element $a+u$ runs over the coset $a+E$, whose
weights are $\{0,2,2,2\}$ if $a\in E$ and $\{1,1,1,3\}$ if $a\in O$. We
distinguish four cases according to the parities of $a$ and $b$.

\smallskip
\emph{Case 1: $a,b\in E$.} Since $c\notin H$ we have $a\ne b$. Choosing $u=a\in
E$ in \eqref{eq:cosetcond} gives $\wt(0)+\wt(a+b)=\wt(a+b)$. But
$a+b\in E\setminus\{0\}$ has weight $2$, contradicting \eqref{eq:cosetcond}.

\smallskip
\emph{Case 2: $a\in E$, $b\in O$.} Choosing $u=a$ gives $\wt(a+b)\ge 3$. As
$a+b\in O$ has weight $1$ or $3$, this forces $a+b=111$, i.e.,\ $b=a+111$.
Conversely, if $b=a+111$ then for every $u\in E$ we have
$b+u=(a+u)+111=\overline{a+u}$, the componentwise complement, whence
\[
  \wt(a+u)+\wt(b+u)=3 \qquad\text{for all } u\in E,
\]
so \eqref{eq:cosetcond} holds and $d(\C)=3$. There are four admissible choices
$a\in E$, but they all determine the same code: if $a,a'\in E$ and
$c=(a,a+111)$, $c'=(a',a'+111)$, then $c'=c+(u_0,u_0)$ with $u_0=a+a'\in E$, so
$c'\in \mathcal H+c$. Taking $a=000$ gives $c=000111$ and
\[
  \mathcal H+c=\{000111,\;011100,\;101010,\;110001\},
\]
hence $\C=\CA$.

\smallskip
\emph{Case 3: $a\in O$, $b\in E$.} Symmetrically $a=b+111$, and taking $b=000$
gives $c=111000$ and
\[
  \mathcal H+c=\{111000,\;100011,\;010101,\;001110\},
\]
hence $\C=\CA^{\perp}$.

\smallskip
\emph{Case 4: $a,b\in O$.} For $u\in E$ both $a+u$ and $b+u$ lie in $O$ and
therefore have weight $1$ or $3$. Hence \eqref{eq:cosetcond} holds for a given
$u$ if and only if $a+u=111$ or $b+u=111$, that is, if and only if
$u\in\{a+111,\;b+111\}$. This set has at most two elements, whereas $|E|=4$, so
\eqref{eq:cosetcond} fails for some $u\in E$.

\smallskip
Cases 2 and 3 exhaust the possibilities, giving exactly the two codes stated.
Both contain $\mathcal H$, and $\CA\cap\CA^{\perp}=\mathcal H$.
\end{proof}

\begin{corollary}\label{cor:rho}
The overturning $\tau$ interchanges the two trigram-coherent $[6,3,3]$ codes:
\[
  \tau(\CA)=\CA^{\perp},\qquad \tau(\CA^{\perp})=\CA .
\]
Consequently the unordered pair $\{\CA,\CA^{\perp}\}$, and hence the tetrad
$\mathcal H=\CA\cap\CA^{\perp}$, is canonically attached to the trigram decomposition of
$\FF_2^6$: it is determined by the parameters $[6,3,3]$ together with the
requirement of trigram coherence, with no further choice involved.
In particular,
$\CA$ is an isodual code which is not self-dual, the permutation realising the isoduality being the overturning $\tau$.
\end{corollary}

\begin{proof}
The first assertion is clear since $\tau$ sends the generators of $\CA$ into $\CA^{\perp}$. The second follows since $\tau$
is an involution. Alternatively, $\tau$ preserves both the parameters
$[6,3,3]$ and the subspace $\Delta$, so it permutes the set of
trigram-coherent codes, which by Theorem~\ref{thm:uniqueness} has exactly two
elements; since $\CA$ is not self-dual, $\tau$ cannot fix $\CA$.
\end{proof}

\begin{remark}\label{rem:count}
An exhaustive computation in Magma/Python confirms Theorem~\ref{thm:uniqueness} and places it
in context. There are exactly $30$ binary $[6,3,3]$ codes in $\FF_2^6$; they
form a single orbit under the coordinate group $S_6$, so all of them are
equivalent to the shortened Hamming code and share the weight enumerator
$x^6+4x^3y^3+3x^2y^4$, and each has a stabiliser of order
$720/30=24$ in $S_6$. Of these $30$ codes, exactly $2$ contain $\mathcal H$, namely
$\CA$ and $\CA^{\perp}$. Equivalence therefore says nothing about the present
question: what singles out $\CA$ and $\CA^{\perp}$ is not the abstract code but
its position relative to the fixed trigram splitting
$\FF_2^6=\FF_2^3\oplus\FF_2^3$, which the group $S_6$ does not preserve.
\end{remark}

\begin{remark}[I-Ching reading]\label{rem:iching}
Proposition~\ref{prop:bound} says that the eight doubled-trigram hexagrams
cannot all be accommodated by a single-error-correcting core: at most four of
them may coexist in a $[6,3,3]$ code, and the four that can are forced to be
\[
  \{\,\text{Earth},\;\text{Fire},\;\text{Wind},\;\text{Lake}\,\},
\]
because these are exactly the doubled trigrams $(u,u)$ with $u$ of even weight.
The complementary tetrad
$\{\text{Heaven},\text{Water},\text{Mountain},\text{Thunder}\}$ is excluded for
a structural reason: each of $001001$, $010010$, $100100$ has Hamming weight
$2$, so no minimum-distance-$3$ code can contain it, while $111111$ lies at
distance $2$ from $\mathcal H$. Thus these four hexagrams are pushed into
the deep-hole coset $M_A$. Thus the splitting of the eight doubled trigrams
into a coded tetrad and a deep-hole tetrad is not a feature of one convenient
example: it is the only splitting the parameters allow.
\end{remark}

\section{Deep holes for $\CA$} \label{deepholes}

Recall that the distance from a vector $x\in\mathbb F_2^6$ to the code
$\CA$ is
\[
        d(x,\CA)=\min\{d(x,c):c\in \CA\}.
\]
The covering radius of $\C_A$ is
\[
        \rho(\CA)=\max_{x\in\mathbb F_2^6} d(x,\CA).
\]
A vector $x$ is called a {\it deep hole} of $\CA$ if
\[
        d(x, \CA)=\rho(\CA).
\]

The following is the case
$r=3$ of the classical fact~\cite[p.533]{MacWilliamsSloane} that shortened Hamming $[2^r-2,2^r-r-2,3]$ codes are nearly perfect, and so quasi-perfect.

\begin{proposition} \label{prop-quasi-perfect}
$\rho(\CA) =2$, that is, $\CA$ is quasi-perfect. The set of deep holes of $\CA$ is
\[
        M_A=\C_A+(e_3+e_6).
\]
\end{proposition}

\begin{proof}
For the code $\C_A$, the first seven cosets
\[
        \C_A,\quad \C_A+e_1,\quad \dots,\quad \C_A+e_6
\]
consist precisely of the vectors at distance $0$ or $1$ from $\C_A$. Since $d(\CA) = 3$, they are disjoint and the coset leader of $\C_A+e_i$ $(1\le i \le 6)$ has weight $1$.
They cover $56$ vectors. The remaining coset is
\[
        M_A=\C_A+(e_3+e_6),
\]
whose coset leader has weight $2$ since all the cosets of weight 1 are given above. Hence every vector in $M_A$ has
distance exactly $2$ from $\C_A$. Therefore
\[
        \rho(\C_A)=2.
\]
Consequently,
\[
        M_A
        =
        \{x\in\mathbb F_2^6 : d(x, \C_A)=2\}
\]
is exactly the set of deep holes of $\C_A$.

In particular, the deep holes are
\[
\boxed{
\begin{aligned}
\{&
001001,\ 001110,\ 010010,\ 010101,\\
&100011,\ 100100,\ 111000,\ 111111
\}.
\end{aligned}
}
\]
\end{proof}

\begin{figure}[t]
    \centering
    \includegraphics[width=1.0\linewidth]{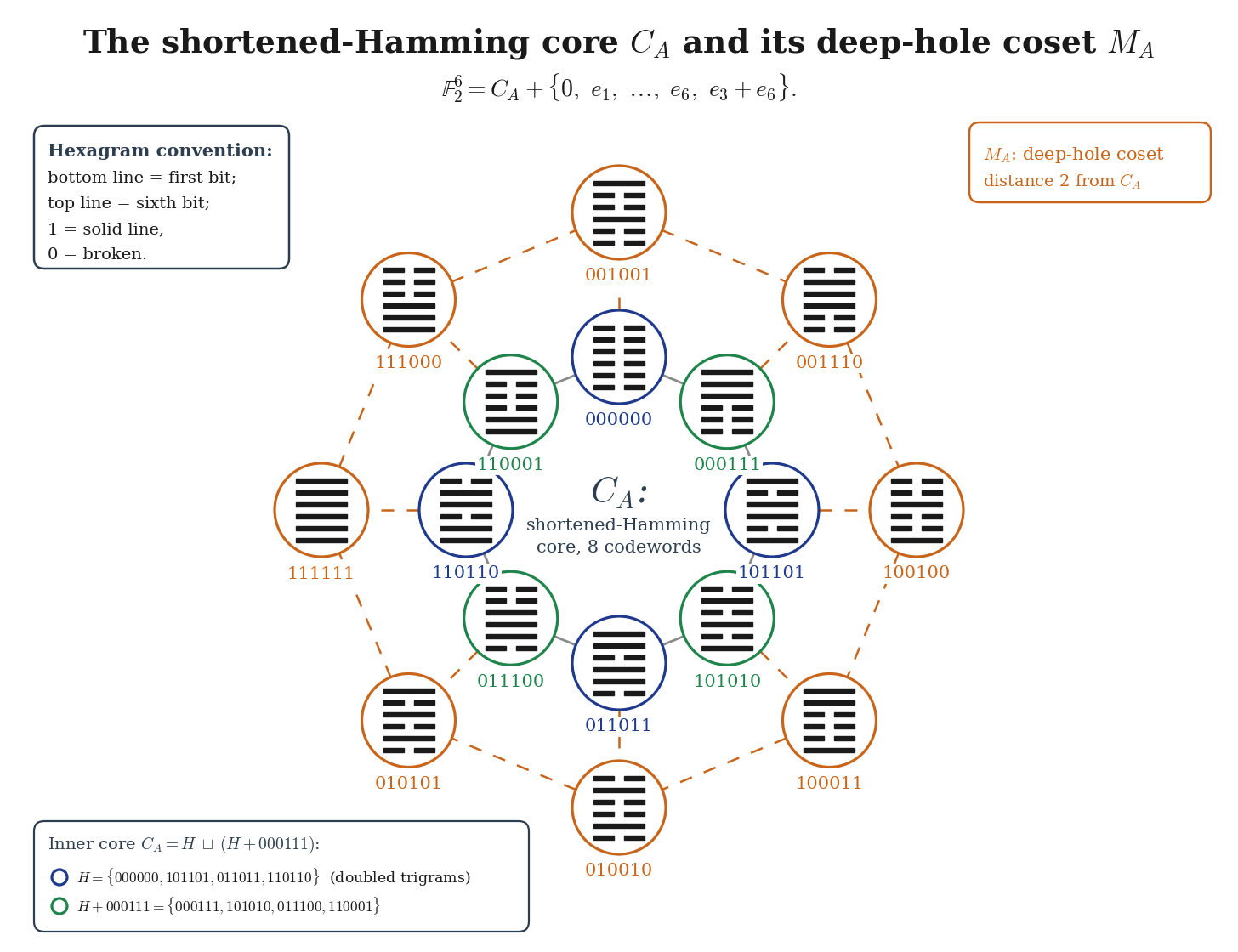}
    \caption{Eight core hexagrams set $\CA$ in the center with the eight orange circled deep holes set $M_A$ outside. Within $\CA$, the four blue circled hexagrams form a $[6,2,4]$ code and the four green circled hexagrams form its coset.}
    \label{fig:hexa_deep-hole}
\end{figure}


The code $\CA$ contains the four doubled-trigram hexagrams
\[
        000000,\qquad 101101,\qquad 011011,\qquad 110110.
\]
These correspond to
\[
        \text{Earth/Earth},\qquad
        \text{Fire/Fire},\qquad
        \text{Wind/Wind},\qquad
        \text{Lake/Lake}.
\]
Thus $\CA$ contains the doubled-trigram tetrad
\[
        \boxed{
        \{\text{Earth},\text{Fire},\text{Wind},\text{Lake}\}.
        }
\]

The missing coset $M_A$ contains the complementary doubled-trigram tetrad
\[
        111111,\qquad 010010,\qquad 001001,\qquad 100100,
\]
corresponding to
\[
        \text{Heaven/Heaven},\qquad
        \text{Water/Water},\qquad
        \text{Mountain/Mountain},\qquad
        \text{Thunder/Thunder}.
\]
Hence $M_A$ contains the complementary tetrad
\[
        \boxed{
        \{\text{Heaven},\text{Water},\text{Mountain},\text{Thunder}\}.
        }
\]

We visualize $\CA$ and its deep-hole coset in Figure~\ref{fig:hexa_deep-hole}. The blue circled hexagrams are generated by $101101, 011011$ to form a diamond while the green circled hexagrams are obtained by adding $000111$ to the blue circled hexagrams to form  a square. Each orange circled hexagram has exactly three hexagrams in blue or green circles at distance exactly 2. It is interesting to note that each orange circled hexagram is complementary to a unique hexagram in blue or green circle because $111111$ is another coset representative of $M_A$.

\begin{proposition}
Let $\widehat{\CA}=\CA\cup M_A$. Then
\[
  \widehat{\CA}
  =\CA+\CA^{\perp}
  =\bigl(\CA \cap \CA^{\perp}\bigr)^{\perp}=\mathcal H^{\perp},
\]
a binary $[6,4,2]$ code with weight enumerator
$1+3z^2+8z^3+3z^4+z^6$ and covering radius $\rho(\widehat{\CA})=1$.
Thus $\CA$ alone does not cover $\FF_2^6$ with radius one --- it
misses exactly the eight deep holes --- while adjoining the deep-hole coset
yields a code that does.
\end{proposition}

\begin{proof}
Since $111000\in\CA^{\perp}\setminus\CA$ and
$111000+001001=110001 \in \CA$, we have
$\CA^{\perp}\subseteq \CA\cup M_A$, so
$\CA+\CA^{\perp}=\CA \cup M_A$; this set has
$16$ elements and is a subspace because $001001 \notin \CA$ and therefore
$\CA \cup (\CA + 001001)= \langle \CA, 001001 \rangle_{\F_2}$.
The identity $(\CA \cap \CA^{\perp})^{\perp}
=\CA^{\perp}+\CA$ is standard, and
$\CA \cap \CA^{\perp}=\mathcal H$ by Theorem~\ref{thm:uniqueness}, so
$\widehat{\CA}=\mathcal H^{\perp}$.  The weight enumerator follows by
inspection of the sixteen codewords. Every $x \in \F_2^6$ lies within distance 1 of
$\widehat{\CA}$ because if $d(x, \CA) \le 1$ this is clear, and otherwise $x \in M_A \subseteq \widehat{\CA}$ by Proposition~\ref{prop-quasi-perfect}, so $d(x, \widehat{\CA})=0$. Hence $\rho(\widehat{\CA}) \le 1$, and $\rho(\widehat{\CA}) \ne 0$ since
$\widehat{\CA} \subsetneq \F_2^6$. Therefore, $\rho(\widehat{\CA}) = 1$, as desired.
\end{proof}

\begin{remark}
Now $\mathcal H=\CA \cap \CA^{\perp}$ is the common self-orthogonal core fixed by both dualization and the
overturning $\tau$; the sole asymmetry between $\CA$ and its dual is the exchange of the single
non-doubled generator $000111\leftrightarrow111000$, itself an overturning pair. In particular, the
isoduality of $\C_A$ is realised by a genuine $I$-Ching symmetry --- reading hexagrams upside down
--- and not merely by an abstract coordinate permutation.
\end{remark}

\section{Codes coherent with a diagonal subspace}
\label{sec:general}

The uniqueness statement of Section~\ref{subsec:uniqueness} concerns a single code of length $6$.
In this section we show that it is the case $m=3$ of a general phenomenon.
Throughout, $m\ge 2$ is an integer and we work in $\FF_2^{2m}$, whose
coordinates are split into a first and a second block of length $m$.

We write $[n,k,d]$ for a binary linear code of length $n$, dimension $k$ and
minimum distance exactly $d$, and $[n,k,\ge d]$ when only the inequality
$d(\C)\ge d$ is asserted. The distinction matters here: the codes we
characterise are those with $d(\C)\ge 3$, and some of them turn out to have
minimum distance $4$ (Corollary~\ref{cor:d4}).

\subsection{The diagonal subspace and a general bound}

\begin{definition}\label{def:diag}
The \emph{doubling map} is the injective linear map
\[
  \iot:\FF_2^{m}\longrightarrow\FF_2^{2m},\qquad \iot(u)=(u,u),
\]
and its image $\Delta_m=\iot(\FF_2^m)$ is the \emph{diagonal subspace} of
$\FF_2^{2m}$. Clearly $\dim\Delta_m=m$ and
\begin{equation}\label{eq:wtdouble}
  \wt(\iot(u))=2\wt(u)\qquad (u\in\FF_2^m).
\end{equation}
\end{definition}

For $m=3$, the diagonal subspace is the doubled-trigram subspace $\Delta$ (Section~\ref{subsec:uniqueness}), and
$\iot$ is the map sending a trigram to the hexagram obtained by repeating it.

The following observation reduces the question ``how much of $\Delta_m$ can a
code of minimum distance $d$ contain?'' to a standard code-parameter problem
one dimension down. Write $k_2(m,\delta)$ for the largest dimension of a binary
linear code of length $m$ and minimum distance at least $\delta$.

\begin{proposition}\label{prop:genbound}
Let $\C\le\FF_2^{2m}$ be a linear code with $d(\C)\ge d$. Then
\[
  T:=\iot^{-1}\!\left(\C\cap\Delta_m\right)
\]
is a binary linear code of length $m$ with $d(T)\ge\lceil d/2\rceil$.
Consequently
\[
  \bigl|\, \C\cap\Delta_m\,\bigr| \;\le\; 2^{\,k_2\left(m,\;\lceil d/2\rceil\right)} .
\]
\end{proposition}

\begin{proof}
Since $\iot$ is an injective linear map and $\C\cap\Delta_m$ is a subspace of
$\Delta_m$, the preimage $T$ is a subspace of $\FF_2^m$ and
$|\C\cap\Delta_m|=|T|$. Let $u\in T$ be nonzero. Then $\iot(u)$ is a nonzero
codeword of $\C$, so $\wt(\iot(u))\ge d$, and by \eqref{eq:wtdouble} this reads
$2\wt(u)\ge d$, i.e.,\ $\wt(u)\ge d/2$. As $\wt(u)$ is an integer,
$\wt(u)\ge\lceil d/2\rceil$. Hence $d(T)\ge\lceil d/2\rceil$ and $\dim T\le
k_2(m,\lceil d/2\rceil)$.
\end{proof}

We now specialise to $d=3$, the case of single-error-correcting codes, where
$\lceil d/2\rceil=2$ and the bound becomes completely explicit. Let
\[
  E_m=\{\,u\in\FF_2^m:\wt(u)\equiv 0 \bmod 2\,\}
\]
denote the even-weight code, an $[m,m-1,2]$ code.

\begin{corollary}\label{cor:d3bound}
Let $\C\le\FF_2^{2m}$ be a linear code with $d(\C)\ge 3$. Then
\[
  \bigl|\, \C\cap\Delta_m\,\bigr|\;\le\;2^{\,m-1},
\]
and equality holds if and only if $\C\cap\Delta_m=\iot(E_m)$.
\end{corollary}

\begin{proof}
By Proposition~\ref{prop:genbound}, $T=\iot^{-1}(\C\cap\Delta_m)$ has
$d(T)\ge 2$, i.e.,\ $T$ contains no vector of weight $1$. A binary code of
length $m$ and dimension $m$ is all of $\FF_2^m$ and contains such a vector, so
$\dim T\le m-1$. If $\dim T=m-1$ then $T=\ker\varphi$ for a nonzero linear
functional $\varphi(x)=a\cdot x$, and $T$ avoids the weight-one vectors
$e_1,\dots,e_m$ if and only if $a\cdot e_i=1$ for every $i$, that is, if and
only if $a=\onev$. Hence $T=E_m$ is the unique possibility, and
$\iot(E_m)=\C\cap\Delta_m$.
\end{proof}

\begin{definition}\label{def:coherentgen}
A binary linear code $\C\le\FF_2^{2m}$ with $d(\C)\ge 3$ is
\emph{diagonal-coherent} if it attains the bound of
Corollary~\ref{cor:d3bound}, that is, if $\iot(E_m)\subseteq \C$. We call $\C$
\emph{extremal} if in addition its dimension is as large as possible.
\end{definition}

\subsection{Existence at optimal redundancy}

Requiring diagonal coherence is a strong constraint: it prescribes a subcode of
dimension $m-1$ inside a code of length $2m$. It is therefore not obvious that
a diagonal-coherent code can also be optimal as an error-correcting code. We
show that it always can, and that the constraint in fact costs nothing.

Recall that for binary codes of length $n$ with $d\ge 3$ the sphere-packing
bound forces the redundancy $n-k$ to satisfy $2^{\,n-k}\ge n+1$, and that this
is attained by shortened Hamming codes. We write
\begin{equation}\label{eq:r0}
  r_0(n)=\lceil\log_2(n+1)\rceil
\end{equation}
for the resulting optimal redundancy, and note for later use that the
minimality of $r_0(n)$ is equivalent to the two-sided estimate
\begin{equation}\label{eq:range}
  2^{\,r_0(n)-1}\;\le\; n\;\le\; 2^{\,r_0(n)}-1 .
\end{equation}

\begin{lemma}\label{lem:dual}
Let $\mathcal D=\iot(E_m)\le\FF_2^{2m}$. Then $\dim \mathcal D=m-1$ and
\[
  \mathcal D^{\perp}=\bigl\{\,(a,\,a+\varepsilon\onev)\;:\;a\in\FF_2^m,\ \varepsilon\in\FF_2\,\bigr\},
  \qquad \dim \mathcal D^{\perp}=m+1 .
\]
\end{lemma}

\begin{proof}
A vector $y=(a,b)$ satisfies $y\cdot\iot(u)=a\cdot u+b\cdot u=(a+b)\cdot u$, so
$y\in \mathcal D^\perp$ if and only if $a+b\in E_m^{\perp}$. The dual of the even-weight
code is the repetition code, $E_m^{\perp}=\{0,\onev\}$, which gives the stated
description; the dimension count is $2m-(m-1)=m+1$.
\end{proof}

\begin{theorem}\label{thm:nocost}
Let $m\ge 2$ and set $n=2m$ and $r=r_0(n)=\lceil\log_2(2m+1)\rceil$. Then there
exists a binary $[\,2m,\ 2m-r,\ \ge 3\,]$ code containing $\iot(E_m)$. In other
words, the maximum dimension of a diagonal-coherent code of length $2m$ equals
the maximum dimension of \emph{any} binary code of length $2m$ with minimum
distance at least $3$: diagonal coherence entails no loss of redundancy.
\end{theorem}

\begin{proof}
Put $\mathcal D=\iot(E_m)$. A code $\C$ of length $2m$ satisfies $\mathcal D\subseteq \C$ and
$d(\C)\ge 3$ if and only if $\mathcal R:=\C^{\perp}$ is contained in $\mathcal D^{\perp}$ and the
columns of a generator matrix of $\mathcal R$ --- equivalently, of a parity-check matrix
of $\C$ --- are pairwise distinct and nonzero.

Let $\mathcal R\le \mathcal D^{\perp}$ have dimension $s$ and let $M$ be an $s\times 2m$ matrix
whose rows form a basis of $\mathcal R$. By Lemma~\ref{lem:dual} the $i$-th row has the
form $(a^{(i)},\,a^{(i)}+\varepsilon_i\onev)$. Set
\[
  e=(\varepsilon_1,\dots,\varepsilon_s)^{\!\top}\in\FF_2^{s},
  \qquad
  c_j=(a^{(1)}_j,\dots,a^{(s)}_j)^{\!\top}\in\FF_2^{s}\quad(1\le j\le m).
\]
Then the columns of $M$ are exactly
\begin{equation}\label{eq:columns}
  c_1,\dots,c_m,\;c_1+e,\dots,c_m+e .
\end{equation}
Let $\pi:\FF_2^{s}\to\FF_2^{s}/\langle e\rangle$ be the quotient map. The
columns \eqref{eq:columns} are pairwise distinct and nonzero if and only if
\begin{equation}\label{eq:cond}
  e\ne 0,\qquad
  \pi(c_1),\dots,\pi(c_m)\ \text{are pairwise distinct and}\ \ne\pi(0).
\end{equation}
Indeed, $e\ne 0$ is needed to separate $c_j$ from $c_j+e$; distinctness of the
$\pi(c_j)$ is exactly the condition that no two of the $2m$ columns coincide;
and $\pi(c_j)\ne\pi(0)$ says $c_j\notin\{0,e\}$, i.e.\ neither $c_j$ nor
$c_j+e$ vanishes.

The quotient $\FF_2^{s}/\langle e\rangle$ has $2^{s-1}-1$ nonzero elements, so
\eqref{eq:cond} can be satisfied if and only if
\begin{equation}\label{eq:sizecond}
  m\le 2^{s-1}-1,\qquad\text{that is,}\qquad 2^{s}\ge 2m+2 .
\end{equation}
Since $2^{s}$ is even and $2m+1$ is odd, the inequality $2^{s}\ge 2m+1$ is
equivalent to $2^{s}\ge 2m+2$. Hence the least $s$ for which \eqref{eq:cond} is
solvable is exactly $s=\lceil\log_2(2m+1)\rceil=r$.

It remains to check that a solution with $s=r$ really produces a code of
dimension $2m-r$, i.e.,\ that $M$ has rank $r$. Choose the $\pi(c_j)$ to be $m$
distinct nonzero elements of $\FF_2^{r}/\langle e\rangle$ and lift them
arbitrarily. Any proper subspace of the $(r-1)$-dimensional quotient contains
at most $2^{r-2}-1$ nonzero elements; from \eqref{eq:range} we have
$2^{r-1}\le 2m$, hence $2^{r-2}-1\le m-1<m$, so the $m$ distinct nonzero
elements $\pi(c_j)$ cannot lie in a proper subspace and therefore span
$\FF_2^{r}/\langle e\rangle$. Together with $e\in\langle$columns$\rangle$ this
shows that the columns of $M$ span $\FF_2^{r}$, so $\operatorname{rank}M=r$ and
$\C=\mathcal R^{\perp}$ has dimension $2m-r$ and $d(\C)\ge 3$. Finally
$\iot(E_m)\subseteq \C$ because $\mathcal R\subseteq \mathcal D^{\perp}$, and equality
$\C\cap\Delta_m=\iot(E_m)$ then holds by Corollary~\ref{cor:d3bound}.
\end{proof}

\subsection{Enumeration}

The proof of Theorem~\ref{thm:nocost} is constructive enough to count the
extremal codes exactly.

\begin{theorem}\label{thm:count}
Let $m\ge2$, $r=\lceil\log_2(2m+1)\rceil$, and let $N(m)$ denote the number of
extremal diagonal-coherent codes, that is, the number of binary
$[\,2m,\ 2m-r,\ \ge3\,]$ codes containing $\iot(E_m)$. Then
\[
  N(m)\;=\;\frac{\displaystyle (2^{r}-1)\cdot 2^{m}\cdot\prod_{i=1}^{m}\bigl(2^{\,r-1}-i\bigr)}
                {\displaystyle \prod_{i=0}^{r-1}\bigl(2^{r}-2^{i}\bigr)} ,
\]
the denominator being $\bigl|\GL_r(\FF_2)\bigr|$.
\end{theorem}

\begin{proof}
By the proof of Theorem~\ref{thm:nocost}, such codes $\C$ correspond bijectively
to subspaces $\mathcal R=\C^{\perp}\le D^{\perp}$ of dimension $r$ satisfying
\eqref{eq:cond}, and a choice of basis of $\mathcal R$ is the same as a tuple
\[
  \alpha=(e;\,c_1,\dots,c_m),\qquad e\in\FF_2^r\setminus\{0\},\ c_j\in\FF_2^r,
\]
subject to \eqref{eq:cond}. The number of such tuples is
\[
  (2^r-1)\cdot\prod_{i=1}^{m}\bigl(2^{\,r-1}-i\bigr)\cdot 2^{m},
\]
since $e$ may be any nonzero vector, the images $\pi(c_1),\dots,\pi(c_m)$ form
an ordered sequence of distinct nonzero elements of a group of order
$2^{r-1}$, and each $\pi(c_j)$ has exactly two preimages.

Two tuples determine the same subspace $\mathcal R$ precisely when they are related by a
change of basis, i.e.,\ by the natural action of $\GL_r(\FF_2)$,
\[
  g\cdot(e;c_1,\dots,c_m)=(ge;\,gc_1,\dots,gc_m).
\]
This action is free: if $g$ fixes $\alpha$ then $g$ fixes $e$ and every $c_j$,
and these vectors span $\FF_2^r$ (as shown at the end of the proof of
Theorem~\ref{thm:nocost}), so $g$ is the identity. Hence the number of
subspaces, and therefore of codes, is the number of tuples divided by
$|\GL_r(\FF_2)|=\prod_{i=0}^{r-1}(2^r-2^i)$.
\end{proof}

\subsection{The minimum distances of the extremal diagonal-coherent codes}

An extremal diagonal-coherent code has $d\ge3$ by definition, but not
necessarily $d=3$. It turns out that the possible values of $d$ are governed by
a fact about \emph{all} binary codes of minimal redundancy, in which the
diagonal plays no role; we prove that first and then specialise.

Recall that a subset $S$ of an abelian group is \emph{sum-free} if
$x+y\notin S$ for all $x,y\in S$ with $x\ne y$.

\begin{lemma}\label{lem:sumfree}
Let $S\subseteq\FF_2^{r}\setminus\{0\}$ be sum-free. Then $|S|\le 2^{\,r-1}$,
and equality holds if and only if $S=\FF_2^r\setminus K$ for some subgroup $K$
of index~$2$; equivalently, $S=\{x:\lambda(x)=1\}$ for a nonzero linear
functional $\lambda$.
\end{lemma}

\begin{proof}
Fix $s\in S$. For $x\in S\setminus\{s\}$ we have $x+s\notin S$ by sum-freeness,
and $x+s\ne 0$ since $x\ne s$. Hence $S$ and $(S\setminus\{s\})+s$ are disjoint
subsets of $\FF_2^r\setminus\{0\}$, giving $|S|+(|S|-1)\le 2^r-1$, that is,
$|S|\le 2^{r-1}$.

Suppose $|S|=2^{r-1}$. Then the above inclusion is an equality, so
$(S\setminus\{s\})+s=\bigl(\FF_2^r\setminus\{0\}\bigr)\setminus S$. Adjoining
the element $0=s+s$ to both sides gives $S+s=\FF_2^r\setminus S$ for
\emph{every} $s\in S$. Consequently $S+s=S+s'$ for all $s,s'\in S$, whence
$S+(s+s')=S$. Let $K=\langle s+s' : s,s'\in S\rangle$. Then $S+K=S$, so $S$ is a
union of cosets of $K$; for a fixed $s_0 \in S$ and for every $s \in S$, $s+s_0 \in K$, thus $S \subseteq s_o + K$ and since $S+K=S$, we have $s_0 + K \subseteq S$, proving $S = s_0 + K$.
So $|S|=|K|=2^{r-1}$ and $K$ has index $2$.

Finally
$0\notin S$ forces $S=\FF_2^r\setminus K$. The converse is immediate: if
$\lambda(x)=\lambda(y)=1$ then $\lambda(x+y)=0$, so $\{x:\lambda(x)=1\}$ is
sum-free of size $2^{r-1}$.
\end{proof}

The following is essentially folklore — it is equivalent to the classical description of maximal sum-free sets in elementary abelian 2-groups — but we include a short self-contained proof, since it is the equality case that we shall need in Corollary~\ref{cor:d4}.


\begin{theorem}\label{thm:maxd}
Let $n\ge4$ and let $r=r_0(n)=\lceil\log_2(n+1)\rceil$ be the minimal
redundancy of a binary code of length $n$ with minimum distance at least $3$.
Then the largest minimum distance of a binary $[\,n,\,n-r\,]$ code is
\[
  d_{\max}(n)=
  \begin{cases}
    4, & n=2^{\,r-1}, \text{ i.e.\ } n \text{ is a power of } 2,\\[2pt]
    3, & \text{otherwise.}
  \end{cases}
\]
Moreover, when $n=2^{\,r-1}$ the codes attaining $d=4$ are exactly the extended
Hamming codes $[\,2^{\,r-1},\,2^{\,r-1}-r,\,4\,]$, unique up to a permutation
of the coordinates.
\end{theorem}

\begin{proof}
Let $\C$ be a binary $[n,n-r]$ code and let $M$ be a parity-check matrix, an
$r\times n$ matrix of rank $r$ with column set $S$. Then $d(\C)\ge 3$ if and
only if the columns are pairwise distinct and nonzero, and, given that,
$d(\C)\ge 4$ if and only if no three distinct columns sum to $0$, i.e.,\ if and
only if $S$ is a sum-free subset of $\FF_2^r\setminus\{0\}$.

Suppose $d(\C)\ge4$. By Lemma~\ref{lem:sumfree}, $n=|S|\le 2^{r-1}$; combined
with the lower estimate $2^{r-1}\le n$ of \eqref{eq:range} this forces
$n=2^{r-1}$. Since \eqref{eq:range} confines $n$ to the interval
$[2^{r-1},2^{r}-1]$, whose only power of $2$ is $2^{r-1}$, the condition
$n=2^{r-1}$ is the same as $n$ being a power of $2$. The equality case of
Lemma~\ref{lem:sumfree} further gives $S=\{x:\lambda(x)=1\}$ for a nonzero
functional $\lambda$; after a change of basis this is the standard
parity-check matrix of the extended Hamming code of length $2^{r-1}$, whence
the uniqueness assertion. This code has $d=4$ and no more: since
$|S|=2^{r-1}\ge4$ we may pick distinct $x,y,z\in S$, and then $x+y+z\in S$
because $\lambda(x+y+z)=1$, while $x+y+z\notin\{x,y,z\}$; these four distinct
columns sum to $0$.

Conversely, if $n=2^{r-1}$ then taking $S=\{x:\lambda(x)=1\}$ as the column set
produces a rank-$r$ matrix, and the corresponding code has $d=4$ by the above.

Finally, for every $n$ in the range \eqref{eq:range} a shortened Hamming code
realises $d=3$, so $d_{\max}(n)=3$ when $n$ is not a power of $2$.
\end{proof}

\begin{remark}\label{rem:hypothesis}
The minimality of $r$ in Theorem~\ref{thm:maxd} cannot be dropped. For
$n=10$ and $r=5$ --- one more than the minimal value $r_0(10)=4$ --- choose as
columns any $10$ of the $16$ vectors of an affine hyperplane
$\{x\in\FF_2^5:\lambda(x)=1\}$. No three of them sum to $0$, so the resulting
code is a $[10,5,4]$ code although $10$ is not a power of $2$. What
Theorem~\ref{thm:maxd} asserts is that this freedom disappears once the
redundancy is pushed down to $r_0(n)$.
\end{remark}

\begin{corollary}\label{cor:d4}
Let $m\ge2$ and put $r=r_0(2m)$, and let $\C$ be an extremal diagonal-coherent
code of length $2m$.
\begin{enumerate}
\item[(i)] If $d(\C)\ge4$, then $2m=2^{\,r-1}$; in particular $m$ is a power of
      $2$, and moreover $d(\C)=4$ and $\C$ is an extended Hamming code with parameters
      $[\,2m,\,2m-r,\,4\,]$.
\item[(ii)] Conversely, if $m$ is a power of $2$, then at least one extremal
      diagonal-coherent code of length $2m$ attains $d=4$. Not all of them need
      do so: for $m=4$ exactly two of the $N(4)=10$ extremal codes have $d=4$,
      the remaining eight having $d=3$ (Proposition~\ref{prop:countd4}).
\item[(iii)] If $m$ is not a power of $2$, then $d(\C)=3$.
\end{enumerate}
\end{corollary}

\begin{proof}
By Definition~\ref{def:coherentgen} and Theorem~\ref{thm:nocost}, the code $\C$
has length $n=2m$, minimum distance $d(\C)\ge3$ and dimension $2m-r$, so its
redundancy is exactly the minimal value $r=r_0(2m)$; and $m\ge2$ gives
$n=2m\ge4$. Theorem~\ref{thm:maxd} therefore applies to $\C$ and yields
$d(\C)\le d_{\max}(2m)$, where $d_{\max}(2m)=4$ if $2m$ is a power of $2$ and
$d_{\max}(2m)=3$ otherwise.

(i) If $d(\C)\ge4$ then $d_{\max}(2m)\ge 4$, so $2m=2^{\,r-1}$ and
$d(\C)=d_{\max}(2m)=4$; the uniqueness clause of Theorem~\ref{thm:maxd}
identifies $\C$ as an extended Hamming code $[\,2^{\,r-1},\,2^{\,r-1}-r,\,4\,]$.
Since $n=2m\ge4$ forces $r\ge3$, we get $m=2^{\,r-2}$, a power of $2$.

(iii) If $m$ is not a power of $2$, neither is $2m$, so $d_{\max}(2m)=3$;
together with $d(\C)\ge3$ this gives $d(\C)=3$.

(ii) Let $m=2^{\,r-2}$ with $r\ge3$, so that $2m=2^{\,r-1}$. Choose a nonzero
linear functional $\lambda$ on $\FF_2^{r}$ and set
$S=\{x\in\FF_2^r:\lambda(x)=1\}$, so that $|S|=2^{\,r-1}=2m$ and $0\notin S$.
Since $|\ker\lambda|=2^{\,r-1}\ge4$, we may choose
$e\in\ker\lambda\setminus\{0\}$. Then $\lambda(x+e)=\lambda(x)$, so
translation by $e$ is a fixed-point-free involution of $S$ and partitions $S$
into $|S|/2=2^{\,r-2}=m$ pairs $\{c,c+e\}$. Pick one representative $c_j$ from
each pair. The data $(e;c_1,\dots,c_m)$ satisfies \eqref{eq:cond}: $e\ne0$ by
choice; $\pi(c_i)=\pi(c_j)$ would mean $c_j\in\{c_i,c_i+e\}$, contradicting
that distinct $c_j$ represent distinct pairs; and $\pi(c_j)\ne\pi(0)$ because
$\lambda(c_j)=1$ while $\lambda(0)=\lambda(e)=0$. By the proof of
Theorem~\ref{thm:nocost} these data produce an extremal diagonal-coherent code
$\C_0$ of length $2m$ whose parity-check columns \eqref{eq:columns} are exactly
the $2m$ elements of $S$. As $\lambda(x)+\lambda(y)=0$ for $x,y\in S$, no two
elements of $S$ sum to a third, i.e.,\ $S$ is sum-free; by the proof of
Theorem~\ref{thm:maxd} this says $d(\C_0)\ge4$, whence $d(\C_0)=4$ by~(i). The
count for $m=4$ is Proposition~\ref{prop:countd4} below.
\end{proof}


\begin{corollary}\label{cor:optimal}
Let $m\ge2$, let $r=r_0(2m)$, and let $\C$ be an extremal diagonal-coherent code
of length $2m$ whose minimum distance is the largest occurring among such
codes, that is, $d(\C)=4$ if $m$ is a power of $2$ and $d(\C)=3$ otherwise. Then
$\C$ is optimal in both of the following senses.
\begin{enumerate}
\item[(i)] \emph{Distance-optimal}: no binary $[\,2m,\,2m-r\,]$ code has
      minimum distance larger than $d(\C)$.
\item[(ii)] \emph{Dimension-optimal}: no binary code of length $2m$ with
      minimum distance at least $d(\C)$ has dimension larger than $2m-r$.
\end{enumerate}
Consequently the starred entries of the minimum-distance column of Table~\ref{tbl:N(m)} record optimal parameters, and for $m=3$ this recovers the distance-optimality
of $\CA$ noted in Section~\ref{sec-shortHamming}.
\end{corollary}

\begin{proof}
By Definition~\ref{def:coherentgen} and Theorem~\ref{thm:nocost}, $\C$ is a
binary $[\,2m,\,2m-r\,]$ code, and $2m\ge4$.

(i) Theorem~\ref{thm:maxd}, applied with $n=2m$, states that the largest
minimum distance of a binary $[\,2m,\,2m-r\,]$ code is $d_{\max}(2m)$, equal to
$4$ when $2m$ is a power of $2$ and to $3$ otherwise. Since $2m$ is a power of
$2$ if and only if $m$ is, Corollary~\ref{cor:d4} says exactly that $\C$ attains
$d(\C)=d_{\max}(2m)$.

(ii) Suppose first $d(\C)=3$, and let $\C'$ be a binary code of length $2m$ with
$d(\C')\ge3$ and redundancy $s=2m-\dim \C'$. The sphere-packing bound gives
$2^{s}\ge 2m+1$, so $s\ge\lceil\log_2(2m+1)\rceil=r$ by \eqref{eq:r0} and
$\dim \C'\le 2m-r$.

Suppose now $d(\C)=4$, so that $2m=2^{\,r-1}$ by Corollary~\ref{cor:d4}, and let
$C'$ be a binary code of length $2m$ with $d(\C')\ge4$ and redundancy $s$. Let
$M'$ be a parity-check matrix of $C'$, of rank $s$, and let
$S'\subseteq\FF_2^{\,s}$ be its set of columns. As in the proof of
Theorem~\ref{thm:maxd}, $d(\C')\ge4$ forces the columns to be pairwise distinct
and nonzero with no three of them summing to $0$; that is, $S'$ is a sum-free
subset of $\FF_2^{\,s}\setminus\{0\}$ of size $2m$. Lemma~\ref{lem:sumfree}
gives $2m\le 2^{\,s-1}$, and $2m=2^{\,r-1}$ therefore yields $s\ge r$, i.e.,
$\dim \C'\le 2m-r$.
\end{proof}

\begin{remark}\label{rem:notoptimal}
When $m$ is a power of $2$, the other extremal diagonal-coherent codes are
\emph{not} distance-optimal: for $m=4$ the eight codes with $d=3$ have
parameters $[8,4,3]$, whereas the two extended Hamming codes show that
parameters $[8,4,4]$ exist. Thus diagonal coherence alone does not force optimality of
the minimum distance; what Corollary~\ref{cor:optimal} says is that the best
minimum distance available under diagonal coherence is already the best
available without it.
\end{remark}

\begin{proposition}\label{prop:countd4}
Let $r\ge3$ and $m=2^{\,r-2}$, so that $2m=2^{\,r-1}$ and $r_0(2m)=r$. Then the
number $N_4(m)$ of extremal diagonal-coherent codes of length $2m$ with
minimum distance $4$ is
\[
  N_4(m)\;=\;\frac{\bigl(2^{r}-1\bigr)\bigl(2^{\,r-1}-1\bigr)\cdot 2^{m}\cdot m!}
                  {\bigl|\GL_r(\FF_2)\bigr|}\,,
\]
all remaining $N(m)-N_4(m)$ extremal codes having minimum distance $3$. In
particular
\[
  N_4(2)=1=N(2),\qquad N_4(4)=2<N(4)=10,\qquad N_4(8)=480<N(8)=205\,920 .
\]
\end{proposition}

\begin{proof}
By the proofs of Theorems~\ref{thm:nocost} and~\ref{thm:count}, the extremal
diagonal-coherent codes of length $2m$ correspond bijectively to the orbits of
the free $\GL_r(\FF_2)$-action on the set of tuples
$\alpha=(e;c_1,\dots,c_m)$ satisfying \eqref{eq:cond}, the code attached to
$\alpha$ having as parity-check columns the set
\[
  S(\alpha)=\{c_1,\dots,c_m,\,c_1+e,\dots,c_m+e\}\subseteq\FF_2^r\setminus\{0\},
  \qquad |S(\alpha)|=2m=2^{\,r-1}.
\]
By the proof of Theorem~\ref{thm:maxd}, the code has $d\ge4$ if and only if
$S(\alpha)$ is sum-free, and since $|S(\alpha)|=2^{\,r-1}$ the equality case of
Lemma~\ref{lem:sumfree} applies: this holds if and only if
$S(\alpha)=\{x:\lambda(x)=1\}$ for some nonzero linear functional $\lambda$ on
$\FF_2^{r}$. By Corollary~\ref{cor:d4}(i) such a code has $d=4$ exactly.

It remains to count the tuples $\alpha$ with $S(\alpha)$ of this form. The
functional $\lambda$ is determined by $S(\alpha)$ and may be chosen in $2^r-1$
ways. Given $\lambda$, both $c_j$ and $c_j+e$ lie in $S(\alpha)$, so
$\lambda(e)=0$; as $e\ne0$ this leaves $2^{\,r-1}-1$ choices of $e$.
Given $\lambda$ and $e$, translation by $e$ partitions $S=\{x:\lambda(x)=1\}$
into $m=2^{\,r-2}$ pairs $\{c,c+e\}$, and a tuple with $S(\alpha)=S$ is precisely
an ordered choice of one representative from each pair: there are $2^{m}$ ways
to choose the representatives and $m!$ ways to order them. Each tuple so
obtained satisfies \eqref{eq:cond}, by the verification carried out in the
proof of Corollary~\ref{cor:d4}(ii). Hence the number of admissible tuples is
$(2^r-1)(2^{\,r-1}-1)\,2^{m}\,m!$.

Finally the $\GL_r(\FF_2)$-action preserves this set of tuples, since
$S(g\cdot\alpha)=g\,S(\alpha)=\{x:\lambda(g^{-1}x)=1\}$, and it is free; dividing
by $|\GL_r(\FF_2)|$ gives the stated formula. The three displayed values are
obtained by substituting $r=3,4,5$.
\end{proof}



\subsection{Normal forms for the extremal diagonal-coherent codes}

\begin{proposition}\label{prop:normalform}
Let $m\ge2$ and $r=r_0(2m)$. For an $(r-1)\times m$ matrix $A$ over $\FF_2$ with
columns $v_1,\dots,v_m$ and a vector $\eps=(\eps_1,\dots,\eps_m)\in\FF_2^{\,m}$
put
\begin{equation}\label{eq:normalform}
  M(A,\eps)\;=\;
  \begin{pmatrix}
     A & A\\[2pt]
     \eps_1\ \cdots\ \eps_m & \eps_1+1\ \cdots\ \eps_m+1
  \end{pmatrix}
  \;\in\;\FF_2^{\,r\times 2m}.
\end{equation}
\begin{enumerate}
\item[(i)] If $v_1,\dots,v_m$ are pairwise distinct and nonzero, then
      $M(A,\eps)$ has rank $r$ and $\C=\ker M(A,\eps)$ is an extremal
      diagonal-coherent code of length $2m$.
\item[(ii)] Conversely, every extremal diagonal-coherent code of length $2m$
      equals $\ker M(A,\eps)$ for some such $A$ and some $\eps$.
\item[(iii)] $d(\C)=4$ if and only if $\{v_1,\dots,v_m\}$ is sum-free in
      $\FF_2^{\,r-1}$, and $d(\C)=3$ otherwise. In particular $d(\C)$ does not
      depend on $\eps$.
\end{enumerate}
\end{proposition}

\begin{proof}
Write $e=(0,\dots,0,1)^{\!\top}\in\FF_2^{\,r}$. The columns of $M(A,\eps)$ are
\[
  c_j=\binom{v_j}{\eps_j}\ (1\le j\le m),
  \qquad
  c_j+e=\binom{v_j}{\eps_j+1},
\]
so they are of the shape \eqref{eq:columns} with this $e$, and condition
\eqref{eq:cond} reads exactly: $v_1,\dots,v_m$ pairwise distinct and nonzero.

(i) The column set contains $c_j$ and $c_j+e$, hence contains $e$ in its span.
Modulo $\langle e\rangle$ the columns give $v_1,\dots,v_m$; a proper subspace of
$\FF_2^{\,r-1}$ has at most $2^{\,r-2}-1$ nonzero elements, while
\eqref{eq:range} gives $2^{\,r-1}\le 2m$, i.e. $m\ge2^{\,r-2}$, so the $v_j$
span $\FF_2^{\,r-1}$. Hence the columns span $\FF_2^{\,r}$ and
$\rank M(A,\eps)=r$. The proof of Theorem~\ref{thm:nocost} now shows that
$\C=\ker M(A,\eps)$ has dimension $2m-r$, satisfies $d(\C)\ge3$ and contains
$\iot(E_m)$.

(ii) Let $\C$ be extremal diagonal-coherent and let $M$ be a parity-check matrix
whose rows form a basis of $\mathcal R=\C^{\perp}\le \mathcal D^{\perp}$, with associated data
$(e';c_1',\dots,c_m')$ as in the proof of Theorem~\ref{thm:nocost}. Since
$e'\ne0$ there is $g\in\GL_r(\FF_2)$ with $ge'=e$. Replacing $M$ by $gM$ does
not change the row space, hence not the code; and the $j$-th column of $gM$ is
$gc_j'=\binom{v_j}{\eps_j}$ while its $(m+j)$-th column is
$gc_j'+ge'=\binom{v_j}{\eps_j+1}$. Thus $gM=M(A,\eps)$ with
$A=(v_1\ \cdots\ v_m)$, and \eqref{eq:cond} makes the $v_j$ distinct and
nonzero.

(iii) As a set the columns of $M(A,\eps)$ are
$\{v_1,\dots,v_m\}\times\FF_2$, which does not depend on $\eps$; call it $S$.
By the proof of Theorem~\ref{thm:maxd}, $d(\C)\ge4$ if and only if $S$ is
sum-free. Let $(x,a),(y,b)\in S$ be distinct. If $x=y$ then $a\ne b$ and
$(x,a)+(y,b)=e\notin S$, because $0\notin\{v_1,\dots,v_m\}$. If $x\ne y$ then
$(x,a)+(y,b)=(x+y,a+b)$ lies in $S$ precisely when $x+y\in\{v_1,\dots,v_m\}$.
Hence $S$ is sum-free if and only if $\{v_1,\dots,v_m\}$ is. Since $d(\C)\le4$
by Theorem~\ref{thm:maxd}, the remaining case gives $d(\C)=3$.
\end{proof}

\begin{remark}\label{rem:normalform-count}
Proposition~\ref{prop:normalform}(iii) reproves Corollary~\ref{cor:d4}: a
sum-free subset of $\FF_2^{\,r-1}$ has at most $2^{\,r-2}$ elements by
Lemma~\ref{lem:sumfree}, while $m\ge2^{\,r-2}$ as above, so $d(\C)=4$ forces
$m=2^{\,r-2}$, that is $2m=2^{\,r-1}$. It also gives a second route to
Theorem~\ref{thm:count}: there are
$\prod_{i=1}^{m}\bigl(2^{\,r-1}-i\bigr)$ choices of $(v_1,\dots,v_m)$ and
$2^{m}$ of $\eps$, and the matrices $M(A,\eps)$ and $M(A',\eps')$ define the
same code exactly when they differ by an element of the stabiliser of $e$ in
$\GL_r(\FF_2)$, a subgroup of index $2^{\,r}-1$.
\end{remark}

\begin{corollary}\label{cor:explicit3}
Let $m\ge3$ and $r=r_0(2m)$, and let $v_j\in\FF_2^{\,r-1}$ be the binary
expansion of $j$ for $1\le j\le m$; this is legitimate because
$m\le 2^{\,r-1}-1$ by \eqref{eq:range}. Then
\[
  M\;=\;
  \begin{pmatrix}
     A & A\\[2pt]
     0\ \cdots\ 0 & 1\ \cdots\ 1
  \end{pmatrix},
  \qquad A=(v_1\ v_2\ \cdots\ v_m),
\]
is the parity-check matrix of an extremal diagonal-coherent
$[\,2m,\,2m-r,\,3\,]$ code. Indeed $v_1+v_2+v_3=0$, so $\{v_1,\dots,v_m\}$ is
not sum-free and Proposition~\ref{prop:normalform}(iii) applies.
\end{corollary}

\begin{remark}\label{rem:m3normal}
For $m=3$ we have $r=3$ and $(v_1,v_2,v_3)$ runs through all of
$\FF_2^{2}\setminus\{0\}$, so
\[
  M(A,\eps)=\begin{pmatrix}
    1&0&1&1&0&1\\ 0&1&1&0&1&1\\ \eps_1&\eps_2&\eps_3&\eps_1{+}1&\eps_2{+}1&\eps_3{+}1
  \end{pmatrix}.
\]
Taking $\eps=(0,0,0)$ makes $M(A,\eps)$ equal to the matrix $G_A$ of Section~\ref{sec-shortHamming}. Here it serves as a \emph{parity-check} matrix, so the resulting code is ker$G_A =\CA^{\perp}$ rather than $\CA$ itself; taking $\eps=(1,1,1)$ replaces the last row by $111000$ and yields $\CA$. The eight values of $\eps$ therefore produce only these two codes, in accordance with $N(3)=2$.
\end{remark}


\subsection{Hulls of the extremal diagonal-coherent codes}

Recall that the \emph{hull} of a linear code is
$\Hull(\C)=\C\cap \C^{\perp}$; it measures how far $C$ is from being self-dual at
one extreme and linear complementary dual at the other. For the code
$\CA$ of Section~\ref{sec-shortHamming} we found
$\Hull(\CA)=\mathcal H=\iot(E_3)$ (Theorem~\ref{thm:uniqueness}). The normal
form of Proposition~\ref{prop:normalform} determines the hull in general.

\begin{theorem}\label{thm:hull}
Let $m\ge2$, $r=r_0(2m)$, and let $\C=\ker M(A,\eps)$ be an extremal
diagonal-coherent code written in the normal form \eqref{eq:normalform}. Let
$\Row(A)\le\FF_2^{\,m}$ be the row space of $A$, an $[\,m,\,r-1\,]$ code, and
put
\[
  \sigma \;=\; v_1+v_2+\cdots+v_m \;\in\;\FF_2^{\,r-1},
\]
the sum of the columns of $A$. Then exactly one of the following holds.
\begin{enumerate}
\item[(i)] If $\sigma\ne0$, then
      $\Hull(\C)=\iot\bigl(\Row(A)\cap E_m\bigr)$ and
      $\dim\Hull(\C)=r-2$.
\item[(ii)] If $\sigma=0$ and $m$ is odd, then $\Hull(\C)=\iot\bigl(\Row(A)\bigr)$
      and $\dim\Hull(\C)=r-1$.
\item[(iii)] If $\sigma=0$ and $m$ is even, then $\Hull(\C)=\C^{\perp}$, that is,
      $\C^{\perp}\subseteq \C$, and $\dim\Hull(\C)=r$.
\end{enumerate}
In cases (i) and (ii) the hull is contained in the diagonal subspace
$\Delta_m$. In all cases $\dim\Hull(\C)\ge r-2\ge1$, so no extremal
diagonal-coherent code is linear complementary dual.
\end{theorem}

\begin{proof}
Write $M=M(A,\eps)$ and let $a_1,\dots,a_{r-1}\in\FF_2^{\,m}$ be the rows of
$A$. Since $\rank M=r$ by Proposition~\ref{prop:normalform}(i), the map
$z\mapsto z^{\!\top}M$ is an isomorphism from $\FF_2^{\,r}$ onto
$\C^{\perp}$, and $z^{\!\top}M$ lies in $\C=\ker M$ if and only if
$M\,(z^{\!\top}M)^{\!\top}=MM^{\!\top}z=0$. Hence it restricts to an isomorphism
$\ker\bigl(MM^{\!\top}\bigr)\to\Hull(\C)$, and
$\dim\Hull(\C)=r-\rank\bigl(MM^{\!\top}\bigr)$.

We compute $MM^{\!\top}$. For $1\le i,j\le r-1$ the corresponding rows of $M$
are $(a_i\mid a_i)$ and $(a_j\mid a_j)$, whose inner product is
$2\,(a_i\cdot a_j)=0$: \emph{the doubled rows are pairwise orthogonal and
self-orthogonal}. Writing $\onev$ for the all-one vector of length $m$ and
$\beta_i=a_i\cdot\onev=\wt(a_i)\bmod2$, the inner product of $(a_i\mid a_i)$
with the last row $(\eps\mid\eps+\onev)$ is
$a_i\cdot\eps+a_i\cdot(\eps+\onev)=\beta_i$, and the last row has inner
product $\wt(\eps)+\wt(\eps+\onev)=m$ with itself. Therefore
\[
  MM^{\!\top}=\begin{pmatrix} 0 & \beta\\ \beta^{\!\top} & m \bmod 2\end{pmatrix},
  \qquad \beta=(\beta_1,\dots,\beta_{r-1})^{\!\top}.
\]
The $i$-th entry of $\sigma=\sum_j v_j$ is $\sum_j (v_j)_i=\wt(a_i)\bmod2$, so
$\beta=\sigma$.

If $\sigma\ne0$, then $MM^{\!\top}$ has rank $2$ and its kernel consists of the
$z=(z',z_r)$ with $z_r=0$ and $\sigma\cdot z'=0$. For such $z$ the
corresponding codeword is $z'^{\!\top}(A\mid A)=(w\mid w)=\iot(w)$ with
$w=z'^{\!\top}A\in\Row(A)$, and
$\sigma\cdot z'=\sum_i z_i'\wt(a_i)=\wt(w)\bmod 2$, so the condition on $z'$
says exactly $w\in E_m$. This gives (i), the dimension being
$(r-1)-1=r-2$ because $\Row(A)\not\subseteq E_m$ when $\sigma\ne0$.

If $\sigma=0$ and $m$ is odd, the only constraint is $z_r=0$, giving
$\Hull(\C)=\iot(\Row(A))$ of dimension $r-1$. If $\sigma=0$ and $m$ is even
then $MM^{\!\top}=0$, so $\Hull(\C)=\C^{\perp}$ has dimension $r$.

Finally $2m\ge4$ forces $r\ge3$, whence $\dim\Hull(\C)\ge1$.
\end{proof}

\begin{corollary}\label{cor:hullsimplex}
Suppose $m=2^{\,r-1}-1$, the case treated in Corollary~\ref{cor:hamming} below. Then the
columns of $A$ are all the nonzero vectors of $\FF_2^{\,r-1}$, so $\Row(A)$ is
the simplex code $\mathrm{Sim}_{r-1}$ with parameters $[\,2^{\,r-1}-1,\ r-1,\ 2^{\,r-2}\,]$; moreover $\sigma=0$ and
$m$ is odd. Hence every extremal diagonal-coherent code of length
$2^{\,r}-2$ has
\[
  \Hull(\C)=\iot\bigl(\mathrm{Sim}_{r-1}\bigr),\qquad \dim\Hull(\C)=r-1 .
\]
For $r=3$ this is $\Hull(\CA)=\iot(E_3)=\mathcal H$, since the simplex code of
length $3$ is the even-weight code $E_3$; for $r=4$ it gives the doubled
$[7,3,4]$ simplex code inside $\FF_2^{14}$.
\end{corollary}

\begin{remark}\label{rem:hulltable}
For $3\le m\le6$ case (i) or (ii) of Theorem~\ref{thm:hull} applies to every
extremal code with $d=3$, and the hull dimension is constant for each $m$ in
Table~\ref{tbl:N(m)}: it equals $2,2,2,2$ for $m=3,\dots,6$, and $3$ for
$m=7$. This is an accident of the small values of $r$ involved. For $m=5$ and
$m=6$, for instance, $\sigma$ is the sum of the $7-m$ omitted nonzero vectors
of $\FF_2^{3}$ and is therefore never $0$. At $m=8$, where $r=5$, both regimes
occur among the codes with $d=3$: choosing $(v_j)$ with $\sigma\ne0$ gives
$\dim\Hull(C)=3$, while a set of eight nonzero vectors of $\FF_2^{4}$ summing
to $0$ but not sum-free --- for example
$\{1,2,3,4,5,6,8,15\}$ written as binary expansions --- gives
$\dim\Hull(\C)=5$ and $\C^{\perp}\subseteq \C$.
\end{remark}

\subsection{Deep holes for an extremal diagonal-coherent code}

Section~\ref{deepholes} described the 64 hexagrams as the codewords of $\CA$,
their radius-one balls, and a single coset of deep holes. That description
generalises to every $m$, and the number of deep holes turns out not to depend
on which extremal code is chosen.

\begin{proposition}\label{prop:deepgen}
Let $m\ge2$, let $r=r_0(2m)$, and let $\C$ be an extremal diagonal-coherent code
of length $2m$. Then $\rho(\C)=2$, so $C$ is quasi-perfect, and the deep holes of
$\C$ are the union of exactly
\[
  2^{\,r}-1-2m
\]
cosets of $\C$, each of size $2^{\,2m-r}$; hence $\C$ has
$\bigl(2^{\,r}-1-2m\bigr)\,2^{\,2m-r}$ deep holes. In particular this number
depends only on $m$, and not on the choice of $\C$ among the $N(m)$ extremal
codes.
\end{proposition}

\begin{proof}
Let $M$ be a parity-check matrix of $\C$ and let
$S\subseteq\FF_2^{\,r}\setminus\{0\}$ be its set of columns, so that $|S|=2m$
and the elements of $S$ are pairwise distinct because $d(\C)\ge3$. For a
syndrome $s\in\FF_2^{\,r}$ the coset of $\C$ with syndrome $s$ has coset leaders
of weight $w(s)$, the least number of elements of $S$ summing to $s$; thus
$\rho(\C)=\max_s w(s)$.

We have $w(s)=0$ only for $s=0$, and $w(s)=1$ exactly for the $2m$ elements of
$S$. Since $2m+1$ is odd and greater than $1$, it is not a power of $2$, so the
minimality of $r=\lceil\log_2(2m+1)\rceil$ gives $2^{\,r}\ge 2m+2$; hence at
least one syndrome has $w(s)\ge2$ and $\rho(\C)\ge2$.

Conversely let $s\notin S\cup\{0\}$ and suppose $S\cap(S+s)=\varnothing$. Then
$S$ and $S+s$ are disjoint subsets of $\FF_2^{\,r}$ of size $2m$ each, so
$4m\le 2^{\,r}$; but \eqref{eq:range} gives $2^{\,r-1}\le 2m$, i.e.
$2^{\,r}\le 4m$. Hence $S\sqcup(S+s)=\FF_2^{\,r}$. As $0\notin S$ we get
$0\in S+s$, that is $s\in S$, contrary to the choice of $s$. Therefore
$S\cap(S+s)\ne\varnothing$: there is $x\in S$ with $x+s\in S$, and
$s=x+(x+s)$ shows $w(s)=2$. Thus $\rho(\C)=2$.

Consequently the $2^{\,r}$ syndromes split into one with $w=0$, exactly $2m$
with $w=1$, and the remaining $2^{\,r}-1-2m$ with $w=2$; the last are precisely
the deep-hole cosets. Each coset of $\C$ has $|\C|=2^{\,2m-r}$ elements.
\end{proof}

\begin{remark}\label{rem:onecoset}
The picture of Section~\ref{deepholes}, in which the deep holes form a
\emph{single} coset, occurs exactly when $2^{\,r}-1-2m=1$, that is
$m=2^{\,r-1}-1$ --- precisely the case treated in
Corollary~\ref{cor:hamming}, where the extremal codes are Hamming codes
shortened at one coordinate. For $r=3$ this is $m=3$, the hexagram case, and
for $r=4$ it is $m=7$. So the clean decomposition of $\FF_2^6$ into the eight
balls around $\CA$ and one deep-hole coset is not typical of general $m$: it
reflects the fact that the $[6,3,3]$ code is a singly shortened Hamming code.
\end{remark}


\begin{table}[ht]
\centering
\begin{center}
\begin{tabular}{c c c c c c}
\hline
$m$ & length, dimension & $r$ & $N(m)$ & minimum distances occurring
    & deep-hole cosets \\
\hline
$2$ & $[4,1]$    & $3$ & $1$       & $4^*$                             & $3$  \\
$3$ & $[6,3]$    & $3$ & $2$       & $3^*$                             & $1$  \\
$4$ & $[8,4]$    & $4$ & $10$      & $3$ ($8$ codes), $4^*$ ($2$ codes) & $7$  \\
$5$ & $[10,6]$   & $4$ & $60$      & $3^*$                             & $5$  \\
$6$ & $[12,8]$   & $4$ & $240$     & $3^*$                             & $3$  \\
$7$ & $[14,10]$  & $4$ & $480$     & $3^*$                             & $1$  \\
$8$ & $[16,11]$  & $5$ & $205\,920$& $3$ ($205\,440$ codes), $4^*$ ($480$ codes)
                                                                     & $15$ \\
\hline
\end{tabular}
\end{center}
\caption{Parameters $[2m,\,2m-r,\,d=3,4]$, the number $N(m)$ of extremal
diagonal-coherent codes, and the number $2^{\,r}-1-2m$ of deep-hole cosets of
each of them (Proposition~\ref{prop:deepgen}); every such coset has
$2^{\,2m-r}$ elements. Here $^*$ means that the minimum distance is optimal by Corollary~\ref{cor:optimal}.}
  \label{tbl:N(m)}
\end{table}

\noindent
The minimum-distance column in Table~\ref{tbl:N(m)} is consistent with Corollary~\ref{cor:d4}, whose multiplicities are given by Proposition~\ref{prop:countd4}: the value
$4$ occurs exactly in the rows $m=2,4,8$, where $m$ is a power of~$2$.
The counts $N(m)$ for $2 \le m\le 8$ and $N_4(m)$ for $m=4,8$ have been confirmed by exhaustive computer search in Magma/Python.

\subsection{Specialisations}

\begin{corollary}\label{cor:hamming}
Suppose $m=2^{\,r-1}-1$, so that $2m=2^{r}-2$ and the bound
\eqref{eq:sizecond} is met with equality. Then every extremal
diagonal-coherent code is a Hamming code of length $2^{r}-1$ shortened at a
single coordinate, namely the coordinate whose parity-check column is the
vector $e$ of \eqref{eq:columns}; in particular it has parameters
$[\,2^r-2,\ 2^r-2-r,\ 3\,]$ with minimum distance exactly $3$, and
\[
  N\bigl(2^{\,r-1}-1\bigr)=\frac{(2^r-1)\cdot 2^{\,2^{r-1}-1}\cdot\bigl(2^{\,r-1}-1\bigr)!}
                                {\bigl|\GL_r(\FF_2)\bigr|} .
\]
\end{corollary}

\begin{proof}
When $m=2^{r-1}-1$ the elements $\pi(c_1),\dots,\pi(c_m)$ exhaust the nonzero
elements of $\FF_2^r/\langle e\rangle$, so the $2m=2^r-2$ columns
\eqref{eq:columns} are exactly the nonzero vectors of $\FF_2^r$ other than $e$.
A parity-check matrix whose columns are all nonzero vectors of $\FF_2^r$ defines
the Hamming code of length $2^r-1$; deleting the column $e$ shortens it. The
minimum distance is $3$ and not $4$ by Corollary~\ref{cor:d4}, since
$m=2^{r-1}-1$ is not a power of $2$ for $r\ge3$.
\end{proof}

\begin{remark}\label{rem:howrestrictive}
It is natural to ask how restrictive diagonal coherence really is: is every
optimal $[\,2m,\,2m-r,\,3\,]$ code, $r=r_0(2m)$, diagonal-coherent? Since
containing $\iot(E_m)$ is not preserved by coordinate permutations, the
question must be read up to equivalence, and in that form it has a clean
answer. Let $S\subseteq\FF_2^{\,r}\setminus\{0\}$ be the set of the $2m$
parity-check columns of such a code $C$. Comparing with the normal form
\eqref{eq:normalform}, we see that $C$ is equivalent to a diagonal-coherent
code if and only if
\begin{equation}\label{eq:transl}
  S+e=S\quad\text{for some } e\ne0 ,
\end{equation}
for \eqref{eq:transl} lets one split $S$ into the $m$ pairs
$\{c,c+e\}$ demanded by \eqref{eq:columns}, and conversely the column set of
$M(A,\eps)$ is invariant under $e=(0,\dots,0,1)^{\!\top}$. Equivalently, the
complement $T=\FF_2^{\,r}\setminus S$, which has $2^{\,r}-2m$ elements and
contains $0$, must be invariant under some nonzero translation.

In the situation of Corollary~\ref{cor:hamming}, where $m=2^{\,r-1}-1$, we have
$T=\{0,x\}$ for a single nonzero $x$, and $T+x=T$ holds automatically. Hence
\emph{every} optimal $[\,2^{\,r}-2,\,2^{\,r}-2-r,\,3\,]$ code is equivalent to a
diagonal-coherent one: for these lengths diagonal coherence is no restriction
at all. This covers the hexagram case $m=3$, and quantifies the observation of
Remark~\ref{rem:count} that what distinguishes $\CA$ is not the abstract code
but its position relative to the trigram splitting.

For $m\le 2^{\,r-1}-2$, however, $|T|\ge4$ and one may choose $T$ with trivial
translation stabiliser, so optimal codes that are equivalent to no
diagonal-coherent code exist. For $m=6$ the obstruction is visible at once:
$T=\{0,x,y,z\}$ is translation-invariant only if $x+y+z=0$, so omitting three
linearly independent columns, say
$S=\FF_2^{4}\setminus\{0000,0001,0010,0100\}$, produces an optimal $[12,8,3]$
code that is not equivalent to any diagonal-coherent code. A computer search
shows that the $[10,6,3]$ codes fall into four equivalence classes of which
only one satisfies \eqref{eq:transl}, and the $[12,8,3]$ codes into two classes
of which only one does; for $m=7$ and $m=15$ there is a single class and
\eqref{eq:transl} holds, in accordance with the previous paragraph.
\end{remark}

\begin{corollary}[the case $m=3$]\label{cor:m3}
For $m=3$ we have $r=3$ and $N(3)=2$. The two extremal diagonal-coherent codes
are the $[6,3,3]$ codes
\[
  \CA=\langle 101101,\,011011,\,000111\rangle
  \quad\text{and}\quad
  \CA^{\perp}=\langle 101101,\,011011,\,111000\rangle ,
\]
they satisfy $\CA\cap\CA^{\perp}=\iot(E_3)$, and they are interchanged by the
coordinate reversal $\tau$. This recovers Theorem~\ref{thm:uniqueness}.
\end{corollary}

\begin{remark}\label{rem:three}
Three conditions that are easy to conflate should be kept apart.
\begin{enumerate}
\item[(i)] \emph{Minimum distance $4$ is attained} exactly when $m$ is a power
      of $2$, by Corollary~\ref{cor:d4}; this is an infinite family
      $m=2,4,8,16,\dots$, the codes being extended Hamming codes.
\item[(ii)] \emph{The dimension is half the length}, i.e.\ $2m-r=m$, exactly
      when $m=r=\lceil\log_2(2m+1)\rceil$, which holds only for $m=3$ and
      $m=4$.
\item[(iii)] \emph{Some extremal code is self-dual}, which requires both (i)
      and (ii) and therefore happens only for $m=4$: of the $N(4)=10$ extremal
      codes exactly two are self-dual, and these are the extended Hamming
      $[8,4,4]$ codes.
\end{enumerate}
Thus at $m=3$ the extremal codes form a dual pair exchanged by an involution of
the coordinates, whereas at $m=4$ the dual pair degenerates to a self-dual
code. It would be of interest to know whether this dichotomy has a uniform
explanation.
\end{remark}

\begin{remark}\label{rem:iching-gen}
For $m=3$ the diagonal subspace $\Delta_3$ consists of the eight
doubled-trigram hexagrams of the I-Ching, and Corollary~\ref{cor:d3bound} says
that at most four of them can lie in a single-error-correcting code, the four
being forced to be those $(u,u)$ with $\wt(u)$ even. The general statement is
therefore that a diagonal-coherent code sees exactly the ``even half'' of the
diagonal, for every $m$.
\end{remark}

\subsection{Open problems}

\begin{remark}\label{rem:closed}
Within Definition~\ref{def:coherentgen} the minimum-distance question is
completely settled and nothing remains open. An extremal diagonal-coherent
code has redundancy exactly $r_0(2m)$, so Theorem~\ref{thm:maxd} applies:
its minimum distance is $3$, or $4$, the latter occurring only when
$m$ is a power of 2, and then only for some of the codes (Proposition~\ref{prop:countd4}).
No
value $d\ge5$ can occur at that redundancy. The genuinely open questions
concern the bound of Proposition~\ref{prop:genbound} for the \emph{other}
values of $d$, where the relevant optimal dimension is not $2m-r_0(2m)$ but
the much smaller optimum associated with that $d$; this is
Problem~\ref{pb:dcoherent}.
\end{remark}

\begin{problem}\label{pb:dcoherent}
The bound of Proposition~\ref{prop:genbound} is a different inequality for each
$d$. Call $\C \le \FF_2^{2m}$ with $d(\C)\ge d$ \emph{$d$-diagonal-coherent} if it
attains that bound, i.e.,\ if
$|\C\cap\Delta_m|=2^{k_2(m,\lceil d/2\rceil)}$, and let $\kappa(2m,d)$ denote
the largest dimension of a binary code of length $2m$ with minimum distance at
least $d$. Is there, for every $m$ and every $d$, a $d$-diagonal-coherent code
of dimension $\kappa(2m,d)$? Theorem~\ref{thm:nocost} answers this
affirmatively for $d=3$, where $k_2(m,2)=m-1$.

The first case beyond $d=3$ is already instructive. For $d=5$ we have
$\lceil d/2\rceil=3$, so the bound involves $k_2(m,3)$, itself a Hamming-type
quantity. Taking $m=7$ and the Hamming $[7,4,3]$ code $\mathcal{H}_3$, the doubled
code $\iot(\mathcal{H}_3)$ has dimension $k_2(7,3)=4$ and weights $6$ and $8$, so
it is $5$-diagonal-coherent, and a computer search produces codes of length
$14$ and minimum distance $5$ containing $\iot(\mathcal{H}_3)$ whose dimension
equals $\kappa(14,5)=6$ by~\cite{Grassl}. Decide whether this ``no loss'' behaviour persists in
general, and characterise the extremal $d$-diagonal-coherent codes.
\end{problem}

\begin{problem}\label{pb:sfold}
For $s\ge2$ let $\iot_s:\FF_2^m\to\FF_2^{sm}$, $u\mapsto(u,u,\dots,u)$, so that
$\wt(\iot_s(u))=s\wt(u)$. The argument of Proposition~\ref{prop:genbound}
yields $d\bigl(\iot_s^{-1}(\C\cap\Delta_m^{(s)})\bigr)\ge\lceil d/s\rceil$ for
any $\C \le \FF_2^{sm}$ with $d(\C)\ge d$. Develop the analogue of
Theorem~\ref{thm:nocost} in this setting, and identify the codes attaining the
bound.
\end{problem}

\begin{problem}\label{pb:qary}
Over $\FF_q$ the diagonal $\Delta_m=\{(u,u)\}$ still satisfies
$\wt(\iot(u))=2\wt(u)$, so Proposition~\ref{prop:genbound} holds verbatim. For
$d=3$ the relevant quantity is again the largest dimension of an
$[m,k,\ge2]_q$ code, namely $m-1$, but the codes attaining it are no longer
unique when $q>2$. Determine the $q$-ary analogues of
Theorems~\ref{thm:nocost} and~\ref{thm:count}, and of
Corollary~\ref{cor:d4}. Note that Theorem~\ref{thm:maxd} is genuinely binary:
over $\FF_q$ the columns of a parity-check matrix are only required to be
pairwise linearly independent, and the sum-free argument of
Lemma~\ref{lem:sumfree} has no immediate analogue.
\end{problem}

\medskip





\end{document}